\documentclass[12pt]{article}
\usepackage{latexsym}
\usepackage{amsfonts,color}
\usepackage{amssymb,amsmath,amsthm}
\usepackage{epsfig}
\usepackage{graphicx}

\newcommand{\beq}{\begin{equation}}
\newcommand{\eeq}{\end{equation}}
\newcommand{\bea}{\begin{eqnarray}}
\newcommand{\eea}{\end{eqnarray}}
\newcommand\LM{{\bf L}}
\newcommand\M{{\bf M}}
\newcommand\N{{\bf N}}
\newcommand\al{{\alpha}}
\newcommand\la{{\lambda}}
\newcommand\om{{\omega}}

\newcommand\bu{{\bf {u}}}
\newcommand\bv{{\bf {v}}}
\newcommand\br{{\bf {r}}}

\newcommand\tq{{\tilde{q}}}
\newcommand\tp{{\tilde{p}}}

\newcommand\si{{\sigma}}

\newcommand\ka{{\kappa}}

\newcommand\dd{\mathrm{d}}
\newcommand\bo{\mathcal{P}}

\def\Xint#1{\mathchoice
   {\XXint\displaystyle\textstyle{#1}}%
   {\XXint\textstyle\scriptstyle{#1}}%
   {\XXint\scriptstyle\scriptscriptstyle{#1}}%
   {\XXint\scriptscriptstyle\scriptscriptstyle{#1}}%
   \!\int}
\def\XXint#1#2#3{{\setbox0=\hbox{$#1{#2#3}{\int}$}
     \vcenter{\hbox{$#2#3$}}\kern-.5\wd0}}

\def\dashint{\Xint-}

\begin{document}
\title{Bilinear recurrences and addition formulae for hyperelliptic sigma 
functions} 

\author{Harry W. Braden\thanks{ 
School of Mathematics, University of Edinburgh,  
James Clerk Maxwell Building, Kings Buildings,  
Mayfield Road, Edinburgh EH9 3JZ, U.K.  
~~E-mail: hwb@maths.ed.ac.uk}, 
Victor Z. Enolskii\thanks{ 
Department of Mathematics, Heriot-Watt University,  
Edinburgh EH14 4AS, U.K.   
~~E-mail:
vze@ma.hw.ac.uk } and 
Andrew N.W. Hone\thanks{Institute of Mathematics and Statistics, 
University of Kent, 
Canterbury CT2 7NF, U.K.
~~E-mail:
A.N.W.Hone@kent.ac.uk}  
} 
\maketitle 

\begin{abstract}
The Somos 4 sequences are a family of sequences %, usually 
%defined over $\mathbb{Q}$, 
satisfying a fourth order bilinear 
recurrence relation. In recent work, one of us has  
proved that the general term in such sequences can be expressed 
in terms of  
the Weierstrass sigma function for an associated 
elliptic curve. 
%After giving a brief summary these results, 
Here we derive the analogous family of sequences associated with an 
hyperelliptic curve of genus two 
defined by the affine model 
$y^2=4x^5+c_4 x^4+ \ldots +c_1 x+c_0$. 
We show that the sequences associated with such  
curves satisfy bilinear recurrences of order 8. 
The proof requires an addition formula which involves 
the genus two Kleinian sigma function 
with its argument shifted by the Abelian image of 
the reduced divisor of a single point 
on the curve. 
%The addition formula is a special case of the generalized 
%Frobenius-Stickelberger formula in genus two 
%considered by one of these authors with 
%Eilbeck and Previato, and treated for hyperelliptic curves of all genera 
%$g$ by Onishi. 
The genus two recurrences are related 
to a B\"{a}cklund transformation (BT) for an integrable 
Hamiltonian system, 
namely the discrete case (ii) H\'{e}non-Heiles system. 
%We briefly discuss recurrences associated with hyperelliptic 
%curves of genus $g\geq 3$, and connections with work of Cantor.   
\end{abstract}

%\newpage
\section{Introduction}

In recent work \cite{honeblms}, one of us has considered fourth 
order quadratic recurrences of the form 
\beq
\tau_{n+2}  
\tau_{n-2} 
=\al \, 
\tau_{n+1} 
\tau_{n-1} 
+\beta  
\, 
(\tau_n)^2,
\label{bil}
\eeq
where $\al$ and $\beta$ are constant parameters. Such recurrences 
arise in the theory of elliptic divisibility sequences
\cite{ward1, ward2, shipsey} and their generalizations, the Somos 4 sequences 
\cite{rob, swart}.
In that context, both the parameters $\al$, $\beta$
and the iterates $\tau_n$
are integers, or more
generally take values in $\mathbb{Q}$ or a Galois extension,
and in that case the sequences have applications in
number theory, as they provide a potential source 
of large prime numbers \cite{eew, ems}. Moreover the  Somos 4 sequences, 
defined by a recurrence of the form (\ref{bil}), provide a simple example 
of the Laurent phenomenon: taking the initial data 
$\tau_0,\tau_1,\tau_2,\tau_3$ and the parameters $\al ,\beta$ as 
variables, all subsequent terms $\tau_n$ for $n\geq 4$ in the sequence  
are Laurent polynomials in these variables. 
%(rather than just rational functions  
Fomin and Zelevinsky have proved 
that this remarkable ``Laurentness'' property is shared by a variety of 
other recurrences in one and more dimensions, with applications in 
combinatorics and commutative algebra (see \cite{fz} and references). 

In \cite{honeblms} the following theorem was proved: 

\vspace{.1in}

\noindent {\bf Theorem 1.} {\it The general solution
of the quadratic recurrence relation
(\ref{bil}) %for $\al\neq 0$ 
takes the form
\beq
\tau_n=
A\,B^n\frac{\si (z_0 +n\ka )}{
\si (\ka )^{n^2}},
\label{form}
\eeq
where $\ka$ and $z_0$ are non-zero
complex numbers, the constants $A$ and $B$ are given by
\beq
A=\frac{\tau_0}{\si (z_0 )}, \qquad
B=\frac{\si (\ka )\si (z_0 )\,\tau_1}
{\si (z_0+\ka )\,\tau_0},
%\quad\frac{\tau_1}{\tau_0},
\label{consts}
\eeq
and $\si$ denotes the Weierstrass sigma function
of an associated elliptic curve 
\beq \label{canon}  
y^2=4x^3-g_2x-g_3. 
\eeq 
The values $\ka$, $z_0$  and the invariants $g_2$, $g_3$ are precisely 
determined from the initial data     
$\tau_0,\tau_1,\tau_2,\tau_3$ and the parameters $\al ,\beta$.}

In the next section we summarize some facts about elliptic divisibility 
sequences, Somos 4 sequences and the details of the above theorem. 
In particular 
we explain how the result of Theorem 1 is connected to the 
second order solvable mapping 
\beq
f_{n+1}=\frac{1}{f_{n-1}f_n}\left( \al +\frac{\beta}{f_n}\right)   
\label{nonli}
\eeq
which  is a degenerate case 
of the type of mapping 
studied by Quispel, Roberts and Thompson \cite{qrt}. 
(See also \cite{rgs} for some recent work on the 
global behaviour of real-valued  solutions of such mappings.) 
%We also point out how this Theorem is related to B\"acklund transformations 
%(BTs) for integrable Hamiltonian systems, within the 
%framework applied  in \cite{us1, us2} and put in an algebro-geometric 
%setting in \cite{kv}. 

\noindent {\bf Remark. } The case $\al =0$, which was 
excluded from the statement of the Theorem in \cite{honeblms}, 
corresponds to $\ka $ being a half period, so that 
$\wp (\ka ) $ is  a branch point of $E$, but then 
the formula for $\tau_n$ is has a rather trivial alternating form: 
$\tau_{2k}=\tau_0 (\tau_2 / \tau_0 )^k \beta^{k(k-1)/2}$,
$\tau_{2k+1}=\tau_1 (\tau_3/\tau_1)^k \beta^{k(k-1)/2}$.
The map (\ref{nonli}) is the autonomous version of the 
discrete Painlev\'e I equation (qdPI) 
\beq
\label{qdp1}
f_{n+1}=\frac{1}{f_{n-1}f_n}\left( \al q^n +\frac{\beta}{f_n}\right),
\eeq
which 
has a continuum limit to the
first Painlev\'e equation \cite{gramani, rgs}. The qdPI 
map (\ref{qdp1}) has 
tau-functions that yield
a sequence of $q$-polynomials \cite{side}, and in the autonomous 
case $q=1$ this map reduces to (\ref{nonli}).
Matsutani has 
constructed some particular solutions of (\ref{nonli}) using elliptic 
functions, and has also considered certain higher order 
recurrences associated with genus two 
hyperelliptic functions \cite{matsutani} (see section 3). 
The above Theorem
guarantees that only elliptic functions are necessary
to specify the general solution of the second order 
map (\ref{nonli}). 
 
The result of Theorem 1 can also be understood via the 
addition formula 
\beq 
\frac{\si (z+\ka )\si (z-\ka )}
{\si (z)^2\si (\ka )^2}
=\wp (\ka )-\wp (z ) 
\label{addweier} 
\eeq 
for the forward and backward shifted 
Weierstrass sigma function in terms of the $\wp$ function 
(see e.g. \cite{ww}). In section 3 this leads us to derive a  higher 
order generalization of the recurrence relation (\ref{bil}) 
by considering a suitable addition formula for the Kleinian hyperelliptic    
sigma function associated to a 
curve of   genus  two. The hyperelliptic Kleinian sigma functions are a natural 
extension of Weierstrass elliptic functions to the case of higher 
genus (see e.g. \cite{baker, bel} and references). 
The addition formula we consider is a special case of the 
generalized  Frobenius-Stickelberger formula in \cite{eep, onishi}, which 
is the exact genus two analogue of (\ref{addweier}). The main result of 
our considerations is to derive an eighth order bilinear 
recurrence whose  terms 
are given by an analogue of the formula (\ref{form}). As a corollary, 
we also derive the solution of a family of sixth order nonlinear  
difference equations in terms of Kleinian $\wp$ functions.     
The fourth section explains how this recurrence is related 
to the BT (integrable discretization) of a Hamiltonian 
system with two degrees of freedom, namely 
the integrable case (ii) H\'enon-Heiles  
system; this BT first appeared in \cite{us1, us2}, 
and was put in an algebro-geometric
setting in \cite{kv}.
The extension 
to higher genus is briefly discussed in our concluding section.  

\section{Elliptic divisibility and Somos 4 sequences} 

\setcounter{equation}{0}

The sequence 
\beq 
\label{eds1} 
0,1,-1,-1,-1,2,1,-3,5,7,\ldots 
\eeq 
is an example of an elliptic divisibility sequence. It is obtained 
from the recurrence 
\beq
\tau_{n+2}=
\frac{
\tau_{n+1}\tau_{n-1}+(\tau_n)^2
}{
\tau_{n-2}
},
%\qquad \tau_0=\tau_1=\tau_2=\tau_3=1,
\label{som4}
\eeq
with initial data taken as  
$$ 
\tau_1=1, \qquad \tau_2=\tau_3=\tau_4=-1.  
$$ 
The sequence can be consistently extended backwards 
for negative $n\in\mathbb{Z}$, to give an antisymmetric 
sequence with $\tau_{-n}=-\tau_n$. 
Remarkably, despite the division by $\tau_{n-2}$ at each 
iteration of (\ref{som4}), the subsequent terms of the sequence are 
all integers, and they satisfy the divisibility property 
\beq 
\label{divis} 
\tau_n  |\tau_m \qquad \mathrm{ whenever}  
\qquad n|m. 
\eeq 

More generally Morgan Ward \cite{ward1, ward2} introduced 
a family of such antisymmetric sequences defined by recurrences 
of the form 
\beq 
\label{edseq} 
\tau_{n+2}\tau_{n-2}=(\tau_2)^2\tau_{n+1}\tau_{n-1}  
-\tau_1\tau_3(\tau_n)^2,  
\eeq 
which are  derived by considering sequences of rational points $nP$ 
on an elliptic curve $E$ over $\mathbb{Q}$. To obtain integer 
sequences of this kind it is required that 
$$ 
\tau_0=0, \quad \tau_1=1, \quad \tau_2,\tau_3,\tau_4\in\mathbb{Z} 
\qquad \mathrm{with}\quad \tau_2|\tau_4. 
$$  
Using  the 
addition law on $E$ and     
considering the multiples 
$nP=(x_n,y_n)$ of a single point $P$, 
%% NOT TRUE IN GENERAL writing the multiples of a single point $P$ as 
%$$nP=(x_n,y_n)=(a_n/\tau_n^2,b_n/\tau_n^3) $$ 
%(where in the last expression the coordinates are 
%% written as fractions in lowest terms), REMOVE THIS 
Ward derived 
the bilinear recurrence (\ref{edseq}) for $\tau_n$, 
with the general term being written in 
terms of the sigma function associated with the curve $E$, as 
\beq 
\label{edsform} 
\tau_n=\frac{\si (n\ka )}{\si (\ka )^{n^2}}.  
\eeq 
(See 
also \cite{recs, shipsey}.) 

Using the addition formula (\ref{addweier}) for the  
Weierstrass sigma function, it is a simple exercise to use  
the formula (\ref{edsform}) in order show that 
the terms of the elliptic divisibility sequence 
satisfy the  Hankel determinant relation  
\beq \label{hankel} 
\tau_{n+m}\tau_{n-m}=\left|\begin{array}{cc} \tau_m\tau_{n-1} & 
\tau_{m-1}\tau_n \\ 
\tau_{m+1}\tau_{n} & 
\tau_{m}\tau_{n+1} \end{array}\right|, 
\eeq 
for all $m,n\in\mathbb{Z}$. 
Starting from the Hankel determinant formula, it is then easy to 
prove by induction that all $\tau_n$ are integers with the 
divisibility property (\ref{divis}).

If we consider the same recurrence (\ref{som4}) but instead take initial 
data $\tau_0=\tau_1=\tau_2=\tau_3=1$, then we find the sequence of 
integers   
\beq \label{somos4} 
1,1,1,1,2,3,7,23,59,314,\ldots, 
\eeq 
known as the Somos 4 sequence (see \cite{rob, sloane}). (In fact this 
Somos sequence is just obtained by selecting the odd index terms of the 
elliptic divisibility sequence (\ref{eds1}), 
up to an alternating sign.) 
More generally, following the terminology of \cite{rob, swart}, we refer 
to any sequence defined by a bilinear recurrence of the 
form (\ref{bil}) as a Somos 4 sequence, while the particular sequence 
above  is denoted Somos (4). It turns out that any such sequence 
is associated to a sequence of points $P_0+nP$ on 
an associated elliptic curve
$E$: this fact was proved by algebraic means in 
the thesis of Swart \cite{swart}, which refers to  
unpublished results established independently by both Nelson Stephens 
and Noam Elkies.
% \cite{swart}. 
In \cite{honeblms}, one of us gave an alternative
complex analytic proof, leading to the construction of the functional 
form (\ref{form}) of the general term, as in Theorem 1 above. 

The  approach taken in \cite{honeblms} was to regard equation 
(\ref{bil}) as the bilinear form of an integrable map, analogous 
to the bilinear equation satisfied by the tau function for a 
soliton equation \cite{hir, sato}, and then solve the initial value 
problem for the bilinear equation with specified initial data 
$\tau_0, \tau_1, \tau_2, \tau_3$ and parameters $\al$, $\beta$. 
The quantity $\tau_n$ may be regarded as being the tau function 
for the second order nonlinear map (\ref{nonli}), to which 
it is related by the substitution 
$$ 
f_n=\frac{\tau_{n+1}\tau_{n-1}}{(\tau_n)^2}. 
$$  
The map (\ref{nonli}) has a first integral, given 
by 
$$ 
J:=J(f_{n-1},f_n)=
f_{n-1}f_n+\al\left(\frac{1}{f_{n-1}}+\frac{1}{f_n}\right) 
+\frac{\beta}{f_{n-1}f_n}=J(f_n,f_{n+1}). 
$$
The algebraic 
formula for $J$ itself implies that the pair $(f_{n-1},f_n)$ lies on 
an elliptic curve for all $n$. 
In fact the general solution of the recurrence 
can be written in terms of the Weierstrass $\wp$ function for 
a curve in the canonical form (\ref{canon}), as 
$$ 
f_n=\wp (\ka )-\wp (z_0+n\ka ). 
$$ 
%in terms of the Weierstrass $\wp$ function for 
%a curve (\ref{canon}) in canonical form. 
The construction of this curve $E$ and the points $P_0, P\in E$ solves 
the initial value problem for (\ref{nonli}) which then yields the 
solution (\ref{form}) for the recurrence (\ref{bil}). 
%and As it standsval es 
It is convenient for us to summarize the results of \cite{honeblms} 
by expressing  
the solution of this initial value problem 
in the form of an algorithm, as follows: 

\noindent {\bf Step 1:} Find the backwards iterate $\tau_{-1}$ from the 
initial data, evaluate the quantities $f_0=\tau_{1}\tau_{-1}/(\tau_0)^2$ 
and $f_1=\tau_{2}\tau_{0}/(\tau_1)^2$, and use these 
to calculate the integral $$J=J(f_0,f_1)=\wp''(\kappa)$$. 

\noindent {\bf Step 2:} Use $J$, $\al$, $\beta$ to 
calculate $$\la =\frac{1}{3\al}\left( 
\frac{J^2}{4}-\beta\right)=\wp (\ka ).$$ 
This gives the point $P=(\la ,\mu )=(\wp(\ka ), \wp'(\ka ))\in E$, 
with $\mu=\pm\sqrt{\al}$. 

\noindent {\bf Step 3:} Construct the invariants $g_2,g_3$ of the 
curve $E$ as in (\ref{canon}), from the formulae 
$$ 
g_2=12\la^2-2J, \qquad g_3=4\la^3-g_2\la-\al . 
$$ 

\noindent {\bf Step 4:} Iterate (\ref{nonli}) backwards to obtain 
$f_{-1}$ from $f_0$ and $f_{1}$. Hence find the point $P_0=(\nu ,\xi )=
(\wp(z_0 ), \wp'(z_0 ))\in E$ from the formulae 
$$ 
\nu =\la -f_0, \qquad \xi = \frac{f_0^2 (f_1-f_{-1})}{\mu}. 
$$

\noindent {\bf Step 5:} Calculate the values 
$\ka , z_0\in\mathbb{C}$ 
from the elliptic integrals
$$
\ka =\pm\int_\infty^P \frac{\dd x}{y}, \qquad
z_0= \pm\int_\infty^{P_0} \frac{\dd x}{y};
$$ 
these should be interpreted as
the points in the Jacobian $\mathrm{Jac}(E)$ corresponding to the
points $P,P_0\in E$ respectively.  Note that  
because of the involution $y\to -y$
these values are only defined by to an overall $\pm $ sign, subject to the
constraint that $\wp'(\ka )\wp' (z_0)=\xi\mu=f_0^2 (f_1-f_{-1})$ as in
Step 4. Once $z_0$ and $\ka$ are obtained then $A$ are $B$ are found 
from the formulae (\ref{consts}).

\noindent {\bf Remarks. } It is useful to note that the coefficients 
$\al , \beta$ in the recurrence are given as elliptic functions 
of $\ka$ by 
\beq 
\al =\wp '(\ka)^2, \qquad 
\beta=\wp '(\ka)^2\Big(\wp (2\ka )-\wp (\ka )\Big). 
\label{abpar} 
\eeq 
The above solution of the initial value problem establishes 
an exact correspondence between two sets of six parameters: 
the parameters  
$g_2,g_3,\ka ,z_0, A,B$  that specify the elliptic curve $E$, the 
two  points 
$P,P_0\in E$, and the prefactors in (\ref{form}); 
and the parameters $\al, \beta ,\tau_0, \tau_1,\tau_2, \tau_3$ 
specifying the constant coefficients and initial data 
for the recurrence (\ref{bil}). 
In order to interpret (\ref{nonli}) as an integrable map, it is 
necessary to further specify a symplectic structure \cite{rag, ves}; 
symplectic coordinates and a Lax pair were given in \cite{honeblms}, 
which make (\ref{nonli})  equivalent to the discrete $g=1$ odd 
Mumford system in \cite{kv}.    

As an example of the above algorithm, we present the results for 
the Somos (4) sequence (\ref{somos4}), with 
$\al =\beta =\tau_0=\tau_1=\tau_2=\tau_3=1$. We find $\tau_{-1}=2$, 
so $f_0=2$, $f_1=1$ gives $J=4$ in Step 1. In Steps 2 and 3 we 
have $\la =1$, set $\mu=1$ and then find $g_2=4$, $g_3=-1$, 
and in Step 4 we obtain $f_{-1}=3/4$ so that $\nu =-1$, $\xi =1$. 
Thus the Somos (4) sequence corresponds to the sequence of 
points $P_0+nP=(-1,1)+n(1,1)$ on the curve 
$$ 
E: \qquad y^2=4x^3-4x+1.$$ 
Finally, evaluating the elliptic integrals and 
sigma functions to 9 decimal places 
using the MAPLE computer algebra package 
(version 8), we find that the curve has real and imaginary half-periods
$\omega_1=1.496729323$ and $\omega_3=1.225694691i$ respectively, while      
$$ 
\ka -2\omega_1 =-\int_1^\infty (4t^3-4t+1)^{-\frac{1}{2}}\,dt
=-1.134273216,
$$ 
$$
z_0-2\omega_3 =\int_{-1}^\infty (4t^3-4t+1)^{-\frac{1}{2}}\,dt
=0.204680500-1.225694691i,
$$ 
which yield the other quantities in (\ref{form}) as  
$$
\si (\ka )= 1.555836426, \qquad
A= \frac{1}{\si (z_0)}=0.112724016-0.824911687i 
$$
$$
B= 
\frac{\si (\ka )\si (z_0 )}
{\si (z_0+\ka )}=0.215971963+0.616028193i.
$$
However, the sequence of arguments of
the sigma function can be written
more succinctly as
$$
z_0+n\ka \equiv (2n-3)\hat{z}_0, \qquad
\hat{z}_0=0.929592715+\omega_3,
$$
so that the iterates of the recurrence correspond to the
sequence of points $(2n-3)\tilde{P}$ on the curve E,
where $\tilde{P}=(0,1)$, $P=2\tilde{P}$. 
The full sequence of points $n\tilde{P}$ 
is associated with the elliptic divisibility sequence (\ref{eds1}).

Elliptic divisibility sequences are currently of considerable 
interest due to the fact that large prime numbers can occur therein 
(i.e. $\tau_p$ may be prime when the index $p$ is prime,   
see  \cite{eew, ems, shipsey}). 
Cantor has considered the division  
polynomials for odd hyperelliptic curves \cite{cantor}, corresponding 
to sequences of divisors $n(P-\infty)$, which also satisfy higher 
order recurrences written in terms of Hankel determinants; 
Matsutani has obtained the functional form of these 
division polynomials in genus two \cite{matsupsi}. In the next section 
we shall derive an eighth order bilinear  recurrence associated 
with the sequence of divisors $\mathcal{D}_0+n(P-\infty)$, where 
$\mathcal{D}_0=(P_1-\infty)+(P_2-\infty )$ is the reduced divisor of 
two points on a 
hyperelliptic curve of genus two.

\section{Addition of one point in genus two} 
% and Ermakov-Pinney}

\setcounter{equation}{0}

Let us  
consider an algebraic curve $X$ of genus two defined by the 
affine model
%$y^2=4x^5+\la_4 x^4+ \ldots +\la_1 x+\la_0$.
\beq
X:=\Big\{ (x,y)\, \Big| \, y^2=f(x)\equiv 4x^{5}+\sum_{j=0}^{4}c_j x^j 
\Big\},
\label{curve}
\eeq 
which realizes the curve as a two-sheeted covering of the 
Riemann sphere with $2g+1$ branch points in the complex plane plus 
a single branch point $\infty$ at infinity. 
The vectors of canonical holomorphic differentials and canonical meromorphic 
(second kind) differentials are denoted 
$$ 
\dd \bu =\left(\begin{array}{c} \frac{\dd x}{y} \\  \\ \frac{x\,\dd x}{y} 
\end{array}\right), \qquad 
\dd \br =\left(\begin{array}{c} 
\frac{(12x^3+2c_4x^2+c_3x)\dd x}{4y} \\ \\ \frac{x^2\,\dd x}{y} 
\end{array}\right) 
$$ 
respectively. If we let $(A_1,A_2;B_1,B_2)$ denote the canonical homology 
basis for the compact Riemann surface corresponding to $X$, with 
non-vanishing intersections $A_j\cdot B_k=\delta_{jk}$, then 
the $2\times 2$ matrices of $A$- and $B$-periods are 
given by  
$$ 
2{\bf \om}=\left(\begin{array}{ccc} \oint_{A_1}\frac{\dd x}{y} && 
\oint_{A_2}\frac{\dd x}{y} \\ && \\  
\oint_{A_1}\frac{x\,\dd x}{y} && \oint_{A_2}\frac{x\,\dd x}{y} 
\end{array}\right), 
\qquad 
2{\bf \om}'=\left(\begin{array}{ccc} \oint_{B_1}\frac{\dd x}{y} && 
\oint_{B_2}\frac{\dd x}{y} \\ && \\  
\oint_{B_1}\frac{x\,\dd x}{y} && \oint_{B_2}\frac{x\,\dd x}{y} 
\end{array}\right).
\qquad 
$$      

The Jacobian of $X$ is the complex torus  
$\mathrm{Jac}(X)=\mathbb{C}^2/\Gamma$, where 
$\Gamma=2{\bf \om}\mathbb{Z}^2\oplus 2{\bf \om}'\mathbb{Z}^2$ 
is the lattice generated by 
the periods of canonical holomorphic differentials. 
The elements $(P_1,P_2)$ of the 
symmetric product $(X)^2$ can be identified with 
degree zero divisors $\mathcal{D}=(P_1-\infty)+(P_2-\infty)$, 
which are mapped to $\mathrm{Jac}(X)$ by the Abel map: 
$$ 
\bu=\int_{\infty}^{P_1}\dd\bu +\int_{\infty}^{P_2}\dd\bu \, \in 
\mathrm{Jac}(X)  
$$ 
(where here we are basing the map at $\infty$).  

The Kleinian sigma function $\si (\bu)$, which is a quasiperiodic 
function of $(u_1,u_2)^T=\bu\in\mathbb{C}^2$,
is  the genus two analogue of the Weierstrass sigma function. The 
Kleinian $\zeta$ and $\wp$ functions are defined by 
$$ 
\zeta_j(\bu )=\frac{\partial \log\si (\bu )}{\partial u_j}, \qquad 
j=1,2, 
$$ 
$$ 
\wp_{jk}(\bu )=-\frac{\partial^2 \log\si (\bu )}{\partial u_j\partial u_k},  
\quad 
\wp_{jkl}(\bu )= 
-\frac{\partial^3 \log\si (\bu )}{\partial u_j\partial u_k\partial u_l},  
\qquad j,k=1,2.
$$  
We refer the 
reader to other works such as \cite{baker, bel} for a detailed 
introduction to hyperelliptic curves, Kleinian functions and 
their definition in terms of Riemann theta functions  
(see also \cite{bes, eep, onishi, onishi2} and references).

Vectors $\bv$ in 
the theta divisor  $\Theta$ in Jac($X$)
can be characterized by the fact that
%\beq \label{divisor} 
$$\si (\bv )=0.$$  
%\eeq 
We wish to take a vector $\bv \in \Theta$ in the
%lowest stratum of the
theta divisor, given by
\beq \label{vdef}
\bv =\int_\infty^{(\la ,\mu )}
\dd \bu;
\eeq
so $\bv$ corresponds to the (reduced) divisor of a single point
$(\la ,\mu )\in X$, with $\mu^2=f(\la )$.

In \cite{bel} (Theorem 4.9) it is proved that the Baker function
$$
\Phi : \mathrm{Jac}(X)\times X\rightarrow \mathbb{C},  
$$
defined for $\bu\in\mathrm{Jac}(X)$ and $(\la ,\mu )\in X$ by 
\beq 
\label{ba} 
\Phi (\bu ; (\la ,\mu ) )=\frac{\si (\bv -\bu )}{\si_2 (\bv )\si (\bu )} \,
\exp \left(\dashint_\infty^{(\la ,\mu )} \dd \br^T\bu \right), \qquad 
\bv=\int_\infty^{(\la ,\mu )}\dd \bu ,  
\eeq 
satisfies the Schr\"{o}dinger equation
\beq \label{schrola4} 
(\partial_2^2 - 2\wp_{22})\Phi=(\la +c_4/4)\Phi.
\eeq 
Note that we have chosen a particular normalization for
the Baker function compared with \cite{bel},
including the denominator $\si_2 (\bv)\neq 0$, 
and the principal value 
symbol $\dashint$ in (\ref{ba}) denotes the fact that the integral 
of the meromorphic differential $\dd \br$  
is regularized at infinity.

Let us define two different Baker-Akhiezer functions 
$\Phi_\pm$ related 
by the hyperelliptic involution, as  
$$ 
\Phi_\pm = \Phi(\bu ; (\la ,\pm \mu )).  
$$ 
Then from the proof of Theorem 4.9 in \cite{bel} 
(restricting to $g=2$) we
have that
\beq
\partial_2 \log \Phi_{\pm}= 
\frac{ \pm \mu +\partial_2\bo (\la ;\bu )}{2\bo (\la ;\bu )}, 
\label{bakers}
\eeq
where
\beq
\label{bolza}
\bo (\la ;\bu ) = \la^2-\wp_{22}(\bu )\la  -\wp_{12}(\bu ) 
%\equiv (\la -x_1)(\la -x_2)  
\eeq
is the Bolza polynomial in genus two \cite{bolza}.
Hence it follows that $\bo =\bo (\la ;\bu )$ satisfies the Ermakov-Pinney
equation with respect to derivatives $\partial_2$ 
in the variable $u_2$, namely  
\beq
\label{ep}
\bo (\partial_2^2\bo )-\frac{1}{2}(\partial_2 \bo )^2+2V\bo^2+\frac{\mu^2}{2},
\quad V:=-2\wp_{22}(\bu )-\la -c_4/4.
\eeq
It is a well known classical result of Ermakov 
(see \cite{hone} for references)
that the general solution of the Ermakov-Pinney equation (\ref{ep})
is just given by a product of two solutions of the Schr\"{o}dinger
equation $(\partial_2^2+V)\psi=0$ with Wronskian $\mu$.

Taking the difference of the $\pm$ equations (\ref{bakers}) we have
$$
\frac{W(\Phi_-, \Phi_+)}{\Phi_+\Phi_-}=\frac{\mu}{\bo},
$$
where
$$
W(\Phi_-, \Phi_+)=\left|\begin{array}{cc}\Phi_- &  \Phi_+ \\
\partial_2\Phi_- & \partial_2\Phi_+\end{array}\right|
$$
is the Wronskian. Clearly this must be independent of
$u_2$, but we claim that in fact this Wronskian has
precisely the value $\mu $, which means that
the product 
$
\Phi_+ \Phi_- = \bo.
$
Rewriting this in terms of the sigma function,
we can state the following result. 

\noindent {\bf Proposition.}
{\it The Kleinian sigma function
for a hyperelliptic curve (\ref{curve}) 
of genus two  satisfies the following  formula
for addition of a single point on the curve:
\beq
\frac{\si (\bu +\bv ) \si (\bu -\bv ) }{\si (\bu )^2 \si_2 (\bv )^2 }
=\bo (\la ;\bu )
\label{addn}
\eeq
In the above, $\bu \in \mathrm{Jac}(X)$ is a generic vector in the Jacobian,
$\bv \in \Theta \subset \mathrm{Jac}(X)$ is the image of the single point
$(\la ,\mu )\in X$ under the Abel map, and $\bo$ is the Bolza polynomial
defined by (\ref{bolza}).
}

\noindent {\bf Proof}:  Starting from Baker's addition
theorem for genus two \cite{baker},
\beq \label{baker2} 
\frac{\si (\bu +\bv ) \si (\bu -\bv ) }{\si (\bu )^2 \si (\bv )^2 }
=\wp_{22}(\bu )\wp_{12}(\bv )-\wp_{12}(\bu ) \wp_{22}(\bv )
+\wp_{11}(\bv) -\wp_{11}(\bu ),
\eeq 
where $\bu$, $\bv$ are generic points in $\mathrm{Jac}(X)$,
and multiplying both sides by $\si (\bv )^2/\si_2 (\bv )^2$,
the result follows by taking the limit $\si (\bv )\to 0$ as
$\bv$ tends to the theta divisor. It is necessary to
use the fact (see e.g. \cite{onishi})
that the $x$ coordinate of the point $(\la ,\mu  )\in X$ is given,
in terms of derivatives of the sigma function evaluated on the
theta divisor $\Theta$, by the expression
$$
\la =-\frac{\si_1 (\bv )}{\si_2 (\bv )}.
$$
This follows from the fact that differentiating $\si (\bv )=0$
with respect to
$\la $ gives, by the chain rule,
$$
\frac{dv_1}{d\la }\partial_1\sigma (\bv ) 
+\frac{dv_2}{d\la }\partial_2\sigma (\bv )
=\frac{\si_1 (\bv )}{\mu } +
\frac{\la\si_2 (\bv )}{\mu }=0,
$$
for $\bv\in\Theta$ given by (\ref{vdef}). 
%$=\int^z_\infty d\bu 
$\Box$ 

\noindent {\bf Remark.} Enolskii 
and Gibbons recently calculated the exact analogue of (\ref{addn}) 
in genus three. 
The addition formula (\ref{addn}) is a special 
case of the generalized Frobenius-Stickelberger addition formula 
in genus two considered in \cite{eep, onishi}. 
Onishi has further 
generalized the Frobenius-Stickelberger formula to hyperelliptic sigma 
functions for all genera \cite{onishi2}, and the special case 
of the  formula corresponding to  
addition of one point has been applied to the problem of construction    
of Wannier functions for quasi-periodic finite-gap potentials in \cite{bes}. 

Cantor has constructed the division polynomials for hyperelliptic curves,
and obtained certain recurrence relations for them in the paper \cite{cantor},
where in particular an eighth order bilinear recurrence is found in
genus two.
Up to a suitable normalization, Matsutani has considered
the exact analytic expression for these division polynomials,
which are  equivalent to the sequence of functions
$$ 
a_n=\frac{\si (n\bv )}{\si_2(\bv )^{n^2}},  
$$
known as as hyperelliptic  psi-functions \cite{matsutani, matsupsi}.
In the following theorem, we present a sequence of tau-functions 
that generalize these psi-functions and yet 
satisfy the same recurrence of Somos 8 type.

\noindent {\bf Theorem 2.}
{\it Define the sequence $\{ \tau_n \, | \, n\in \mathbb{Z} \}$  
by 
\beq 
\label{seqg2} 
\tau_n=A\, B^n \frac{\si (\bu +n\bv)}{\si_2(\bv )^{n^2}}, 
\eeq 
where $\bu\in\mathrm{Jac}(X)$ is a generic vector in the Jacobian of the 
genus two curve (\ref{curve}),
$\bv \in \Theta \subset \mathrm{Jac}(X)$ is the image of the single point
$(\la ,\mu)\in X$ under the Abel map, $\si$ denotes the Kleinian sigma 
function of the curve, and $A,B$ are arbitrary constants. Then the terms 
of the sequence satisfy a bilinear recurrence of order 8, given by 
\beq 
\label{recg2} 
\tau_{n+4}\tau_{n-4}=\sum_{j=0}^3 \al_j\,\tau_{n+j}\tau_{n-j},  
\eeq      
where the coefficients $\al_j$ (independent of $n$) are given by 
\beq \label{alcoef}  
\al_1=\frac{\si(6\bv )\si(3\bv )^2}
{\si (4\bv ) \si (2\bv )^2\si_2 (\bv )^{30}},  
\qquad 
\al_3= \frac{\si(3\bv )\si (5\bv)}
{\si (2\bv )\si (4\bv )\si_2(\bv )^{14}}, 
%%%%% minus sign removed here: should be overall plus here!%%%%%
\eeq 
\beq \label{alcoef2} 
\al_2 =\frac{\si (4\bv )^2}{\si_2 (\bv )^{24} \si (2\bv )^2}\left( 1- 
\frac{\si (3\bv )^3\si (5\bv )}{\si (4\bv )^3\si (2 \bv )} \right), 
%\eeq 
 \quad 
%\beq \label{alcoef3} 
\al_0=-\si (6\bv )/(\si (2\bv )\,\si_2(\bv )^{32}).  
\eeq 

} 

\noindent{\bf Proof. } %The formulae for 
%the coefficients $\al_j$ are determined uniquely by 
Substituting 
the expression (\ref{seqg2}) 
into (\ref{recg2}) and using Baker's formula (\ref{baker2}) 
together with the result of the Proposition yields 
an expression of the form 
$$ 
C_{0}(\bv )+ C_{11}(\bv )\wp_{11}(\bu )+ 
C_{12}(\bv )\wp_{12}(\bu )+C_{22}(\bv )\wp_{22}(\bu )=0.    
$$ 
The three functions $\wp_{jk}(\bu )$, $j,k=1,2$
on $\mathrm{Jac}(X)$ 
are not linearly dependent 
(although they do satisfy a nonlinear relation \cite{bel, eep}, 
giving the Kummer surface in $\mathbb{CP}^3$). 
Therefore each of the coefficients $C_0(\bv )$, 
$C_{11}(\bv ),C_{12}(\bv ), C_{22}(\bv )$ must vanish, which 
leads to a linear system for the $\al_j$ as functions of $\bv $. 
This 
determines the above formulae for
the coefficients $\al_j$ uniquely  
in terms of $\sigma$, its
derivative $\sigma_2$, and the $\wp_{jk}$ evaluated
at various multiples of $\bv$. The terms involving $\wp_{jk}$ can be
removed by making use of the addition formulae (\ref{baker2})
and (\ref{addn}) to yield the expressions (\ref{alcoef}) and
(\ref{alcoef2}) in terms of $\si$ and $\si_2$ alone.  
%\beq \label{oldalcoef3}
%\al_0=-\la^2\al_1 -\frac{\si (2\bv )^2}{\si_2(\bv )^8}\,\wp_{11}(2\bv )\al_2
%- \frac{\si (3\bv )^2}{\si_2(\bv )^{18}}\,\wp_{11}(3\bv )\al_3 +
% \frac{\si (4\bv )^2}{\si_2(\bv )^{32}}\,\wp_{11}(4\bv ),
%\eeq
Taking the
limit $\bu\to 0$, these are equivalent to 
Matsutani's expressions
for the
coefficients in the eighth order bilinear
recurrrence for the psi-function (see formula (3.13) in
\cite{matsutani}). 
$\Box$ 
 
\noindent {\bf Corollary.} {\it The sequence of Bolza polynomials 
\beq 
\label{boseq} 
f_n=\bo (\la ;\bu +n\bv  )=\la^2-\wp_{22}(\bu +n\bv )\la - 
\wp_{12}(\bu +n\bv ),  
\eeq 
for $\bu\in\mathrm{Jac}(X)$ and $\bv\in\Theta$, 
satisfies the sixth order nonlinear difference equation 
\beq 
\label{diffg2} 
f_n^4\prod_{k=1}^3(f_{n+k}f_{n-k})^{4-k}=\al_0+ 
\sum_{j=1}^3\al_jf_n^j\prod_{k=1}^{j-1}(f_{n+k}f_{n-k})^{j-k}, 
\eeq 
with the coefficients $\al_j$ as given in equations (\ref{alcoef}), 
and (\ref{alcoef2}).}  

\noindent {\bf Proof of Corollary.} Upon setting  
$$ 
f_n=\frac{\tau_{n+1}\tau_{n-1}}{(\tau_n)^2} 
$$ 
with $\tau_n$ as in (\ref{seqg2}), and using the addition 
formula (\ref{addn}), the result is an immediate consequence 
of Theorem 2. $\Box$   

\noindent {\bf Remarks. } Taking $\bu$ as the Abelian image 
of $(P_1, P_2)\in (X)^2$, the sequence  
(\ref{boseq}) corresponds to the linear flow 
$\bu +n\bv$ in the Jacobian, or equivalently the sequence of divisors 
$\mathcal{D}_n\sim\mathcal{D}_0+n(P-\infty )\sim 
\mathcal{D}_0-n(\hat{P}-\infty )$ with 
$\mathcal{D}_0\sim (P_1-\infty )+(P_2-\infty )$ and $\hat{P}=(\la ,-\mu)$, 
the image of $P=(\la ,\mu)$ under the hyperelliptic involution. The Bolza 
polynomial leads to the solution of 
the Jacobi inversion problem for the curve 
(\ref{curve}) (see Theorem 2.2 in \cite{bel}, and section 4 below), 
so that in particular if 
$\mathcal{D}_n\sim (x_1(n),y_1(n))+(x_2(n),y_2(n))-2\infty$, then we have 
$$ 
f_n=\bo (\la ;\bu +n\bv )=(\la -x_1(n))(\la -x_2(n)), 
\quad 
y_j(n)=-\partial_2 \bo (\la ;\bu +n\bv )|_{\la =x_j}, \quad j=1,2. 
$$ 
Cantor's results in \cite{cantor} concern  
the sequence of reduced divisors $\mathcal{D}_n\sim n(P-\infty)$,
corresponding to the multiples of a single point on an odd
hyperelliptic curve of genus $g$.
In particular for $g=2$ he obtains a bilinear recurrence of order
8, which is the degenerate case $\bu =0$ ($P_1\to\infty$,
$P_2\to\infty$) of our construction, while the analytic derivation 
in that case appears in the work of Matsutani \cite{matsutani, matsupsi}.
The sixth order difference equation 
(\ref{diffg2}) appears as equation (3.15) in \cite{matsutani},
where the special solutions with $\bu =0$ are also presented.  

The sigma functions of genus $g$ odd hyperelliptic curves, 
given by $y^2=f(x)$ with $f$ a polynomial of odd degree  
$2g+1$, are 
known to be tau functions of the Korteweg--deVries  
(KdV) hierarchy 
of partial differential equations (see \cite{bel} for 
instance). 
It is also known that when the curve 
degenerates completely to $y^2=4x^{2g+1}$, the corresponding 
sigma function  degenerates to a polynomial (see \cite{bel2, onishi2}), 
which gives a rational solution of KdV in terms of a Schur function (see 
\cite{am} and chapter 14 in \cite{kac}).   
It is instructive to consider the case when
the curve (\ref{curve}) for $g=2$ degenerates to a 
singular rational curve:
$$
y^2=4x^5.
$$
In that case the Kleinian sigma function degenerates to
the Schur function
\beq
\label{schur}
\si (\bu )=u_1-\frac{u_2^3}{3}, %\qquad \bu =(u_1,u_2)^T,
\eeq
which is the tau function of the three-pole rational solution
of the KdV equation
$$
4\partial_1 V=\partial_2^3 V+6V\partial_2V, \qquad
V=-2\wp_{22}(\bu)=2\partial_2^2\log\si (\bu ) .
$$
The theta divisor consists of vectors of the
form
$$
\bv =\left(\begin{array}{c} \int_\infty^{(\la ,\mu )} \frac{\dd x}{\sqrt{4x^5}} \\
\\
 \int_\infty^{(\la ,\mu )} \frac{x\,\dd x}{\sqrt{4x^5}}\end{array}\right)
=\left(\begin{array}{c} -\frac{1}{3}\la ^{-3/2} \\  \\ -\la^{-1/2}
\end{array}\right) \equiv 
\left(\begin{array}{c} \gamma^3/3 \\ \gamma \end{array}\right), \qquad
\gamma\in \mathbb{C},
$$
satisfying $\si (\bv )=0$. It is trivial to check that
the Schur function satisfies the addition 
formula (\ref{addn}).
Defining $\tau_n$ in terms of the Schur function (\ref{schur}) 
by (\ref{seqg2}), it is easy to verify that this gives 
a particular solution of the eighth order recurrence
(\ref{recg2}) with 
$$ \al_0=-\frac{35}{\gamma^{64}}, \quad
 \al_1=\frac{56}{\gamma^{60}}, \quad
 \al_2=-\frac{28}{\gamma^{48}}, \quad
 \al_3=\frac{8}{\gamma^{28}}. \quad
$$

\section{BT for the case (ii) H\'enon-Heiles system}  

\setcounter{equation}{0}

The integrable case (ii) H\'enon-Heiles system is a system of two degrees 
of freedom defined by the natural Hamiltonian 
\beq 
h_1=\frac{1}{2}(p_1^2+p_2^2)+q_2^3+\frac{1}{2}q_2q_1^2-\frac{1}{2}aq_1^2+cq_2 
-\frac{m^2}{2q_1^2}. 
\label{hham} 
\eeq 
The coordinates $q_j$ and momenta $p_j$ are canonically conjugate, 
and Hamilton's equations 
\beq 
\frac{\dd q_j}{\dd t}=\{ h_1,q_j\}, 
\qquad 
\frac{\dd p_j}{\dd t}=\{ h_1,p_j\}, 
\qquad 
j=1,2 
 \label{hameq} 
\eeq 
are equivalent to the ordinary differential equation for 
travelling wave solutions of the fifth order flow in 
the KdV hierarchy \cite{fordy}.  
The equations of motion (\ref{hameq}) can be written in the 
form of a Lax equation 
$$ 
\frac{\dd \LM}{\dd t}=[\N , \LM ], 
$$ 
where the Lax matrix $\LM$ is  
\beq 
\label{lax} 
\LM (x; q_j,p_j )=\left(\begin{array}{cc} 
\frac{1}{2}p_2 x-\frac{1}{8} p_1q_1 
& \mathcal{B}(x;q_j,p_j) 
%2x^3 +(-q_2+4a)x^2 + \frac{1}{8}\left(q_1^2+4q_2^2-aq_2+16a^2+4c\right)x 
%+\frac{1}{8}\left(p_1^2-\frac{m^2}{q_1^2}\right) 
\\ 
2x^2+(q_2+2a)x-\frac{1}{8}q_1^2 & 
-\frac{1}{2}p_2 x+\frac{1}{8} p_1q_1 \end{array}\right),  
 \eeq 
with 
$$\mathcal{B}(x;q_j,p_j) = 2x^3 +(-q_2+4a)x^2 + 
\frac{1}{8}\left(q_1^2+4q_2^2-aq_2+16a^2+4c\right)x
+\frac{1}{8}\left(p_1^2-\frac{m^2}{q_1^2}\right).$$  
(Note that we have $\LM\to -\frac{(u-a)}{8}\,\LM$, $u\to x+a$ 
compared with reference \cite{us1}.) 
The Lax equation  is the compatibility condition for the 
linear system 
\beq \label{linear} 
\LM \Psi = y\Psi, \qquad \frac{\dd \Psi }{\dd t}=\N \Psi, 
\qquad \N =\left(\begin{array}{cc} 0 & x+a -q_2 \\ 
1 & 0 \end{array}\right). \eeq   
The genus two 
spectral curve is of the precise form (\ref{curve}), namely   
$$ 
\mathrm{det} \, (\LM -y {\bf 1} )=y^2- 
4x^5-12ax^4-(c+12a^2)x^3-\frac{1}{2}h_1x^2- 
\frac{1}{2}h_2x-\frac{m^2}{64}=0, 
$$ 
with 
$$ 
h_2=\frac{1}{4}(q_2+2a)p_1^2-\frac{1}{4}q_1p_1p_2-\frac{1}{32}q_1^4 
- \frac{1}{8}\left(q_2^2-2aq_2 +c+4a^2 \right)q_1^2 
-\frac{m^2(q_2+2a)}{4q_1^2}  
$$ 
being the second independent integral, in involution with $h_1$ 
i.e. $\{h_1,h_2\}=0$. The integral $h_2$ generates a second 
commuting flow 
$$ 
\frac{\dd q_j}{\dd s}=\{ h_2,q_j\},
\qquad
\frac{\dd p_j}{\dd s}=\{ h_2,p_j\},
\qquad
j=1,2.  
$$   

Up to a shift of origin, the time variables $s, t$ can be 
identified with the coordinates $u_1,u_2$ respectively 
on $\mathrm{Jac}(X)$. Using the results of 
Theorem 2.2 in \cite{bel},  
the solution of the H\'enon-Heiles system can be reduced 
to the Jacobi inversion problem 
$$ 
\left(\begin{array}{c} s \\ t\end{array}\right)\equiv 
\bu =\int_\infty^{(x_1,y_1)}\dd\bu + 
\int_\infty^{(x_2,y_2)}\dd\bu , 
$$ 
where the separation coordinates $x_1,x_2$ 
are found from the $(2,1)$ entry in the Lax matrix 
(\ref{lax}), given as a multiple of the Bolza 
polynomial by 
$$ 
2x^2+ (q_2+2a) x- \frac{1}{8} q_1^2 = 
2(x-x_1)(x-x_2)=2 \bo (x;\bu ).  
$$ 
The separation variables $(x_j,y_j)$, $j=1,2$ 
correspond to the reduced divisor $\mathcal{D}_0= 
(x_1,y_1)+(x_2,y_2)-2\infty$, and they are related to the 
Kleinian functions by 
$$ 
\wp_{12}(\bu )=-x_1x_2, \quad 
\wp_{22}(\bu )=x_1+x_2, \quad 
\wp_{221}(\bu )=\frac{x_1y_2-x_2y_1}{x_1-x_2}, 
\quad 
\wp_{222}(\bu )= \frac{y_1-y_2}{x_1-x_2} 
$$ 
(see e.g. \cite{eep}).    
Thus the connection with the Bolza polynomial immediately leads to the 
solution of the H\'enon-Heiles system  
in terms of Kleinian $\wp$ functions, which is 
$$ 
q_1^2=16\wp_{12}(\bu ), \quad q_2=-2\wp_{22}(\bu )-2a, 
\quad q_1p_1=8\wp_{221}(\bu ), \quad p_2=-2\wp_{222}(\bu ). 
$$  

It is shown  in \cite{us1} that the case (ii) H\'enon-Heiles system  
has a B\"acklund transformation (BT) 
with parameter $\la$, which is a symplectic map with 
generating function $F(q_1,\tq_1;\la )$ such 
that $dF=\sum_{j=1,2}p_j\dd q_j-\tp_j\dd \tq_j$. The explicit 
form of the generating function is  
$$ 
F= Z+\frac{m}{2}\log\left(\frac{Z-m}{Z+m}\right) 
+\frac{16}{5}Y^5+4(q_2+\tq_2)Y^3+ 
\Big(2(q_2+\tq_2)^2-2q_2\tq_2+\frac{1}{2}(q_1^2+\tq_1^2)+2c\Big)Y, 
$$   
where $Z(q_j,\tq_j)$ and $Y(q_j ,\tq_j)$ are defined by 
$$
Z^2=m^2+\la q_1^2\tq_1^2, \qquad 
Y^2=\la +a-\frac{1}{2}\left(q_2+\tq_2 \right). 
$$ 
The BT can be realized as a similarity transformation on the 
Lax pair (discrete Lax equation) 
 \beq \label{dlax} 
\tilde{\LM}\M=\LM \M ,
\eeq
where $\tilde{\LM} =\LM (x; \tq_j,\tp_j )$ and  
\beq 
\label{mmat} 
\M=\left(\begin{array}{cc} -Y & Y^2 +x-\la \\ 
1 & -Y \end{array} \right) 
\eeq 
is the elementary Darboux matrix (see \cite{shabat}). 
Clearly from (\ref{dlax}) the BT preserves the spectrum of 
$\LM$, and so maps solutions to solutions. 

In fact, the BT was constructed in \cite{us1} by making use of 
the formulae for the Darboux transformation of the Schr\"odinger 
equation, since the components of $\Psi$ in the linear system (\ref{linear})  
are given by $\Psi=(\psi, \partial_2\psi )^T$ with 
$ (\partial_2^2+q_2)\psi=(x+a)\psi$. Then the quantity $Y$ 
appearing in the Darboux matrix can be given explicitly in terms of 
the Baker function $\Phi_+$ defined in (\ref{ba})  as  
$$ 
Y=\partial_2\log\Phi_+\equiv \partial_2\log \Phi(\bu ;\la ) 
=\zeta_2(\bu -\bv )-\zeta_2(\bu)+ 
\dashint_{\infty}^{(\la ,\mu )}\dd r_2, 
$$ 
and by a simple calculation using (\ref{bakers}) we also have 
$$ 
Y=\frac{ \mu +\partial_2\bo (\la ;\bu )}{2\bo (\la ;\bu )}= 
\frac{ - \mu +\partial_2\bo (\la ;\bu -\bv )}{2\bo (\la ;\bu -\bv )}, 
$$  
which is equivalent to equation (2.9) in \cite{us1}. 

It follows from  general results in 
\cite{kv} that the BT defined by the Darboux matrix (\ref{mmat})       
just gives a shift on the Jacobian, and here we see that 
$$ 
\tq_1^2=16\wp_{12}(\bu -\bv ), \quad \tq_2=-2\wp_{22}(\bu -\bv)-2a,
\quad \tq_1\tp_1=8\wp_{221}(\bu -\bv ), \quad \tp_2=-2\wp_{222}(\bu -\bv).
$$
In terms of divisors, we have the equivalence 
$$\tilde{\mathcal{D}}\sim  
\mathcal{D}_0+(\la ,-\mu )-\infty 
\sim (\tilde{x}_1,\tilde{y}_1) + 
(\tilde{x}_2,\tilde{y}_2)-2\infty\sim \mathcal{D}_{-1}. 
$$ 
Similarly, applying the same BT but using the Baker function 
$\Phi_-$ corresponds to adding the point $P=(\la , \mu)$ 
instead of $\hat{P}=(\la , -\mu)$, which gives the divisor 
$\mathcal{D}_1\sim \mathcal{D}_0+(\la ,\mu )-\infty$. 
In   
other words, the $n$th term in the 
sequence of Bolza polynomials (\ref{boseq}), 
with argument $\bu +n\bv\in \mathrm{Jac}(X)$,  
is the result of $n$ applications of 
the H\'enon-Heiles BT, adding the same point $P$ each time. 
We should point out that the H\'enon-Heiles BT given here 
and in \cite{us1} is 
the same as the $g=2$ discrete odd Mumford system in \cite{kv},   
with a particular choice of symplectic structure. 

\section{Conclusions} 

\setcounter{equation}{0}

We have constructed an eighth order bilinear recurrence relation 
(\ref{recg2}) whose $n$th term is expressed in terms of the 
Kleinian sigma function of an odd hyperelliptic curve (\ref{curve}) 
of genus two. Moreover, we have explained how this is connected 
to the BT for the case (ii) integrable H\'{e}non-Heiles system  
of \cite{us1, us2}: %(see also \cite{kv, fed2, suris} for more on 
%integrable maps): 
each shift $n\to n+1$ in the recurrence corresponds to 
the shift on the Jacobian induced by the BT. 

The question remains of whether we can solve the initial value 
problem for an eighth order recurrence of the form (\ref{recg2}). 
In fact, the expression   
(\ref{seqg2}) depends on at most 10 parameters: the five coefficients 
$c_j$, $j=0, \ldots, 4$ that specify the curve plus the three 
points $P_1,P_2,P$ on the curve to specify the 
sequence of divisors $\mathcal{D}_n\sim (P_1-\infty)+
(P_2-\infty)+n(P-\infty )$, and the two prefactors $A,B$. 
(Actually, the constant $c_4$ can be removed
from the start by making a shift in $x$.)
The general eighth order recurrence 
(\ref{recg2}) would have four parameters $\al_j$, $j=0,1,2,3$ 
in general position and 8 initial data $\tau_j$, $j=0,\ldots,7$.   
In fact, by counting arguments, 
we expect that the general solution of a bilinear recurrence 
consisting of $N+2$ bilinear terms,  
$$ 
\tau_{n+N+1 }\tau_{n-N-1} =\sum_{j=0}^{N}\al_{j}\,
\tau_{n+j}\tau_{n-j},  
$$ 
which is of order $2N+2$, should correspond to 
a sequence of divisors $\mathcal{D}_n\sim 
\sum_{j=1}^N(P_j-\infty)+n(P-\infty )$ on an hyperelliptic 
curve of 
genus $N$, since this solution should depend on $3N+3$ parameters. 
This expectation agrees with the form of vector addition formulae 
for Riemann theta functions obtained by Buchstaber and 
Krichever \cite{bk}, which have precisely $N+2$ terms in genus $N$. 
However, preliminary calculations suggest the need to consider 
a different model for these  
curves, without a branch point at infinity. For example, 
sigma functions for curves of the form 
$$ 
y^2=4x^6+\sum_{j=0}^5c_jx^j 
$$    
should give solutions 
to the Somos 6 recurrence 
$$ 
\tau_{n+3}\tau_{n-3}=\sum_{j=0}^2\al_j\,\tau_{n+j}\tau_{n-j}; 
$$ 
this should correspond to the BT for the $g=2$ even Mumford system 
in \cite{kv}, which has a spectral curve of this type. This agrees 
with recent results of  
van der Poorten \cite{poorten}, which show that  
a certain class of Somos 6 sequences 
arise from the continued fraction expansion of the square root of 
a sextic.  

\noindent {\bf Note added in proof. } 
In a private communication, David Cantor has 
shown one of us that by extending the results of \cite{cantor} 
to divisor sequences $\mathcal{D}_0+n(P-\infty )$ it is 
possible to prove that the sequence of tau functions $\tau_n$ 
given by (\ref{seqg2}) in Theorem 2 satisfies a family of 
{\it trilinear} 
recurrence relations given in terms of a Hankel 
type determinant, that is  
\beq 
\label{trili} 
a_2^2\, a_m\, \tau_n\tau_{n+m}\tau_{n-m}=\left| 
\begin{array}{ccc} 
a_m\tau_{n-2} & a_{m+1}\tau_{n-1} & a_{m+2}\tau_{n} \\ 
a_{m-1}\tau_{n-1} & a_{m}\tau_{n} & a_{m+1}\tau_{n+1} \\ 
a_{m-2}\tau_{n} & a_{m-1}\tau_{n+1} & a_{m}\tau_{n+2} \end{array} 
\right|  
\eeq 
for all $m,n\in\mathbb{Z}$, where 
$$ 
a_m=\frac{\si (m\bv )}{\si_2(\bv )^{m^2}}. 
$$ 
The above formula is the genus two analogue of 
a family of higher recurrences satisfied by 
Somos 4 sequences \cite{swartvdp} 
(cf. also Morgan Ward's identity (\ref{hankel}) for elliptic divisibility 
sequences). 
Up to a suitable choice of normalization, 
these $a_m$ are the genus two division 
polynomials derived by Cantor \cite{cantor}, 
whose functional form has been 
specified precisely in the works \cite{matsutani, matsupsi} 
of Matsutani, where they are referred to as hyperellliptic psi-functions. 
In the limit $\bu \to 0$, when $\tau_n\to a_n$, 
the equation (\ref{trili}) reduces to a Hankel 
formula for these psi-functions obtained by algebraic means 
in \cite{cantor}, which is also presented as formula (3.10) in 
\cite{matsutani}, and rederived using analytic techniques in 
\cite{matsupsi}. Setting $m=3$ in (\ref{trili}) yields a 
sixth order trilinear recurrence for $\tau_n$, while setting $m=4$ 
yields an eighth order trilinear recurrence. As was first pointed 
out to us by Cantor, the bilinear recurrence (\ref{recg2}) 
above then follows as a consequence of these two trilinear identities, 
upon eliminating between the two of them. 
Interestingly, Buchstaber
and Leykin have very recently derived a trilinear differential addition
formula for genus two sigma functions - see Theorem 5.6 in \cite{bule}.

%$\al_0=-\si (6\bv )/(\si (2\bv )\,\si_2(\bv )^{32})$. 

\noindent {\bf Acknowledgments.} AH is grateful to the University 
of Kent 
for supporting
the project
{\it Algebraic curves and functional equations
in mathematical physics} with a Colyer-Fergusson Award, 
which funded his visits to Edinburgh. AH also  
thanks the following people: 
Graham Everest for introducing him to the arithmetic of 
quadratic recurrences; Christine Swart for 
providing a copy of her thesis \cite{swart}; 
and Vadim Kuznetsov, Pol Vanhaecke 
and Orlando Ragnisco for many enlightening conversations about BTs. 
Finally, VZE and AH are very grateful to Shigeki Matsutani 
for providing copies of his papers, and for making some 
very insightful comments about our work.

\end{document}